\documentclass[a4paper,12pt]{amsart}

\usepackage[latin1]{inputenc}
\usepackage[line,matrix,frame,curve,poly]{xy}
\usepackage{color}
\usepackage{url}
\usepackage{mathtools}

\def\Dfn#1{{\sf #1}}
\def\maj{\operatorname{maj}}
\def\imaj{\operatorname{imaj}}
\def\inv{\operatorname{inv}}

\def\Des{\operatorname{Des}}
\def\des{\operatorname{des}}
\def\Asc{\operatorname{Asc}}
\def\asc{\operatorname{asc}}

\def\Set{\operatorname{Set}}

\def\Cat{\operatorname{Cat}}
\def\Aqt{\mathcal{A}_n(q,t)}
\def\Atq{\mathcal{A}_n(t,q)}
\def\Aqq{\mathcal{A}_n(q,q^{-1})}

\def\area{\operatorname{area}}

\def\i{{\bf i}}
\def\iDes{\operatorname{iDes}}
\def\iAsc{\operatorname{iAsc}}

\def\Sn{\mathcal{S}_n}
\def\Sk{\mathcal{S}_k}
\def\Dn{\mathcal{D}_n}

\newtheorem{theorem}{Theorem}[section]
\newtheorem{definition-theorem}[theorem]{Definition-Theorem}
\newtheorem{proposition}[theorem]{Proposition}
\newtheorem{lemma}[theorem]{Lemma}
\newtheorem{corollary}[theorem]{Corollary}

\theoremstyle{definition}
\newtheorem{definition}[theorem]{Definition}

\newtheorem{guess}[theorem]{Guess}
\newtheorem{remark}[theorem]{Remark}

\newtheorem{example}[theorem]{Example}

\definecolor{grau}{rgb}{.5 , .5 , .5}
\definecolor{dunkelgrau}{rgb}{.35 , .35 , .35}
\definecolor{schwarz}{rgb}{0 , 0 , 0}

\begin{document}
	\title{On bijections between 231-avoiding permutations and Dyck paths}
	\author{Christian Stump}
	\address{Fakult\"at f\"ur Mathematik, Universit\"at Wien, Nordbergstra{\ss}e 15, A-1090 Vienna, Austria}
	\email{christian.stump@univie.ac.at}
	\urladdr{http://homepage.univie.ac.at/christian.stump/}
	\date{\today}
	\thanks{Research supported by the Austrian Science Foundation FWF, grant P17563-N13 \lq\lq Macdonald polynomials and q-hypergeometric Series\rq\rq}
	
	\begin{abstract}
		We construct a bijection between $231$-avoiding permutations and Dyck paths that sends the sum of the \emph{major index} and the \emph{inverse major index} of a $231$-avoiding permutation to the \emph{major index} of the corresponding Dyck path. Furthermore, we relate this bijection to others and exhibit a bistatistic on $231$-avoiding permutations which is related to the $q,t$-Catalan numbers.
	\end{abstract}

	\maketitle
	
	\section{Introduction}
		The $q,t$-Catalan numbers $\Cat_n(q,t)$ were defined by A.~Garsia and M.~Haiman in \cite{garsiahaiman} as symmetric rational function in $q$ and $t$. They reduce for $q=t=1$ to the well-known \Dfn{Catalan numbers} $\Cat_n := \frac{1}{n+1}\binom{2n}{n}$; more generally, they reduce for $t=q^{-1}$, up to a power of $q$, to the $q$\Dfn{-Catalan numbers} considered by P.A.~MacMahon in \cite{macmahon},
		$$q^{\binom{n}{2}} \Cat_n(q,q^{-1}) = \frac{1}{[n+1]_q}\begin{bmatrix} 2n \\ n \end{bmatrix}_q.$$
		Here, $[k]_q := \frac{1-q^k}{1-q}$ is the usual $q$-extension of the integer $k$, $[k]_q!:= [1]_q[2]_q \dots [k]_q$ is the $q$-factorial of $k$ and $\left[\begin{smallmatrix} k \\ \ell \end{smallmatrix}\right]_q:=\frac{[k]_q!}{[\ell]_q![k-\ell]_q!}$ is the $q$-binomial coefficient. In \cite{garsiahaglund}, A.~Garsia and J.~Haglund proved that $\Cat_n(q,t)$ is in fact a polynomial with non-negative integer coefficients by providing a combinatorial description of $\Cat_n(q,t)$ in terms of the two statistics \emph{area} and \emph{bounce} on \emph{Dyck paths of semilength} $n$. So far, no bijective proof is known for the fact that $\Cat_n(q,t)$ is symmetric in $q$ and $t$.
		
		\bigskip
		
		In this article, we construct a bijection between $231$-avoiding permutations and Dyck paths which sends the sum of the \emph{major index} and the \emph{inverse major index} on $231$-avoiding permutations to the \emph{major index} of the corresponding Dyck path. Using this bijection, we exhibit a polynomial $\Aqt$ which is related to $\Cat_n(q,t)$ by
		$$\Aqq = \Cat_n(q,q^{-1}),$$
		and show bijectively that $\Aqt$ is symmetric in $q$ and $t$.
		
		Finally, we connect the proposed bijection to a bijection between $132$-avoiding permutations and Dyck paths introduced in \cite{krattenthaler2} by C. Krattenthaler and to a bijection from $231$-avoiding permutations and Dyck paths introduced in \cite{bandlowkillpatrick} by J.~Bandlow and K.~Killpatrick. Using these connections, we obtain several generating function identities as corollaries of the main theorem.
		
		\bigskip
		
		The article is organized as follows: In Section~\ref{sec2}, we define pattern-avoiding permutations and Dyck paths and state some of their basic properties. In Section~\ref{sec3}, we construct the announced bijection between $231$-avoiding permutations and Dyck paths. In Section~\ref{sec4} we describe the bistatistic $(\maj,\binom{n}{2}-\imaj)$ on $231$-avoiding permutations, and we prove the symmetry of the corresponding generating function. In Section~\ref{sec5} we show how the bijection can be described in terms of earlier bijections, and in Section~\ref{sec6} we describe how the bistatistic $(\maj,\binom{n}{2}-\imaj)$ could be connected to the $q,t$-Catalan numbers $\Cat_n(q,t)$.
		
	\section{Basic definitions} \label{sec2}
		
		We denote by $\Sn$ the \Dfn{symmetric group} of all \Dfn{permutations} of the set $[n]:=\{1,\dots,n\}$. Usually, we write a permutation $\sigma \in \Sn$ in \Dfn{one-line notation} as the word $\sigma = [\sigma_1,\dots,\sigma_n]$ of length $n$ where $\sigma_i:=\sigma(i)$ (often, we suppress the commas or even the brackets around). A \Dfn{subword of} $\sigma$ is a subsequence $[\sigma_{i_1},\dots,\sigma_{i_k}]$ of $\sigma$ with $i_1 < \dots < i_k$.
		\begin{definition}
			Let $\sigma \in \Sn$ and $\tau \in \Sk$. $\sigma$ avoids the pattern $\tau$ and is called $\tau$\Dfn{-avoiding} if $\sigma$ does not contain a subword of length $k$ having the same relative order as $\tau$.
		In particular, $\sigma$ is called
			\begin{itemize}
				\item \Dfn{$231$-avoiding} if $\sigma$ does {\bf not} contain a subword $[\sigma_i,\sigma_j,\sigma_k]$ with $\sigma_k < \sigma_i < \sigma_j$,
				\item \Dfn{$312$-avoiding} if $\sigma$ does {\bf not} contain a subword $[\sigma_i,\sigma_j,\sigma_k]$ with $\sigma_j < \sigma_k < \sigma_i$.
			\end{itemize}
		\end{definition}
		For $\tau \in \Sk$, we denote the set of all $\sigma \in \Sn$ which are $\tau$-avoiding by $\Sn(\tau)$.
		\begin{remark}
			In \cite{knuth}, D.E.~Knuth proved that for any $\tau \in \mathcal{S}_3$ the number of $\tau$-avoiding permutations is equal to the $n$-th Catalan number $\Cat_n$.
		\end{remark}
	 	An integer $1 \leq i < n$ is called a \Dfn{descent of} $\sigma$ if $\sigma_i > \sigma_{i+1}$, otherwise $i$ is called an \Dfn{ascent}. As we count $n$ as an ascent, this differs from the usual definition. By $\Des(\sigma)$ (respectively $\Asc(\sigma)$) we denote the set of all descents (respectively ascents) of $\sigma$, and we set $\des(\sigma) := |\Des(\sigma)|$ and $\asc(\sigma) := |\Asc(\sigma)|$. As we will often use the descent set and sometimes the ascent set of the inverse of a permutation, we will denote these sets by $\iDes(\sigma)$ and by $\iAsc(\sigma)$ respectively.
	 	
		We furthermore need the following two elementary involutions on $\Sn$: define $\rho$ to be the involution sending $[\sigma_1,\dots,\sigma_n]$ to $[\sigma_n,\dots,\sigma_1]$ and $\sigma \mapsto \hat{\sigma}$ to be the involution sending a permutation to its inverse. We can easily describe the descent set and the inverse descent set of the images of those involutions: sending a permutation to its inverse interchanges $\Des$ and $\iDes$, and for $\rho$ we have
		\begin{align*}
			\Des(\rho(\sigma)) 	&= [n-1] \setminus \{n-i : i \in \Des(\sigma)\}, \\ 
			\iDes(\rho(\sigma)) &= [n-1] \setminus \iDes(\sigma).
		\end{align*}
		Furthermore, define the \Dfn{major index of} a permutation $\sigma$ by
		$$\maj(\sigma) := \sum_{i \in \Des(\sigma)}{i},$$
		and define its \Dfn{inverse major index} to be the major index of its inverse, $\imaj(\sigma) := \maj(\hat{\sigma})$.
		The inverse major index is discussed in detail by D.~Foata and G.-N.~Han in \cite[Chapters 11, 12]{foatahan}.
		
\begin{figure}
			\centering

\setlength{\unitlength}{1pt}

\begin{picture}(140,140)(0,0)
  \linethickness{.25\unitlength}
  \color{grau}
  \put(0  ,0  ){\line(0,1){120}}
 	\put(20 ,20 ){\line(0,1){100}}
 	\put(40 ,40 ){\line(0,1){80 }}
 	\put(60 ,60 ){\line(0,1){60 }}
 	\put(80 ,80 ){\line(0,1){40 }}
 	\put(100,100){\line(0,1){20 }}
 	
 	\put(0,20 ){\line(1,0){20 }}
 	\put(0,40 ){\line(1,0){40 }}
 	\put(0,60 ){\line(1,0){60 }}
 	\put(0,80 ){\line(1,0){80 }}
 	\put(0,100){\line(1,0){100}}
 	\put(0,120){\line(1,0){120}}
 	
 	\put(0,0){\line(1,1){120}}
 	
 	\color{schwarz}
  \linethickness{1.25\unitlength}
  
 	\put(0  ,0  ){\line(0,1){20.625}}
 	\put(20 ,20 ){\line(0,1){20.625}}
 	\put(20 ,40 ){\line(0,1){20.625}}
 	\put(40 ,60 ){\line(0,1){20.625}}
 	\put(80 ,80 ){\line(0,1){20.625}}
 	\put(100,100){\line(0,1){20.625}}
 	
 	\put(0  ,20 ){\line(1,0){20.625}}
 	\put(20 ,60 ){\line(1,0){20.625}}
 	\put(40 ,80 ){\line(1,0){20.625}}
 	\put(60 ,80 ){\line(1,0){20.625}}
 	\put(80 ,100){\line(1,0){20.625}}
 	\put(100,120){\line(1,0){20.625}}

 	\put(20,20){\circle*{5}}
 	\put(40,60){\circle*{5}}
 	\put(80,80){\circle*{5}}
 	\put(100,100){\circle*{5}}

\end{picture}

			\caption{The Dyck path $01 \hspace{5pt} 001 \hspace{5pt} 011 \hspace{5pt} 01 \hspace{5pt} 01 \in \mathcal{D}_6$.}
			\label{fig3}
		\end{figure}
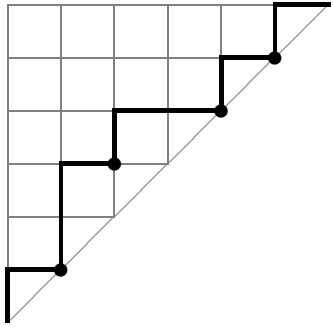

		Next, we define \Dfn{Dyck paths} and state the properties of Dyck paths we need:
		\begin{definition}
			A \Dfn{Dyck path of semilength} $n$ is a sequence of $2n$ letters containing $n$ $0$'s and $n$ $1$'s such that every prefix contains at least as many $0$'s as $1$'s. We denote the set of all Dyck paths of semilength $n$ by $\Dn$.
		\end{definition}
		A natural way of thinking of a Dyck path is to identify it with the lattice path consisting of north and east steps starting at $(0,0)$ and ending at $(n,n)$, where a north step is encoded by a $0$ and an east step by a $1$. The property that any prefix contains at least as many $0$'s as $1$'s is equivalent to the property that a Dyck path always stays weakly above the diagonal $x=y$, see Figure~\ref{fig3} for an example of a Dyck path of semilength $6$.
			
		It is well-known that --~as for $\tau$-avoiding permutations with $\tau \in \mathcal{S}_3$~-- the number of Dyck paths of semilength $n$ is given by $\Cat_n$, see e.g. \cite{stanley2}.
	
		For $D \in \Dn$, an integer $i$ with $1 \leq i < 2n$ is called \Dfn{descent} of $D$, if $D_i = 1$ and $D_{i+1} = 0$, and let $\Des(D)$ denote the set of all descents of $D$. In analogy to the definition on permutations, the \Dfn{major index} of $D$ is defined by
		$$\maj(D) := \sum_{i \in \Des(D)}{i}.$$
		A descent can also be described as a \Dfn{valley} in the lattice path, i.e., an east step followed by a north step. The \Dfn{coordinates} of a descent or valley of a Dyck path $D$ are the coordinates of the lattice point right after the corresponding east step. Furthermore, we define the sets $\Set_X(D)$ and $\Set_Y(D)$ to be the set of $x$-coordinates and the set of y-coordinates of the descents or valleys of $D$. In \cite[Section 3]{callan}, D.~Callan called $\Set_X(D)$ and $\Set_Y(D)$ the \lq\lq ascent-descent code\rq\rq\ of $D$.
		
		The following observation is fundamental.
		\begin{proposition}\label{weneedit3}
			A Dyck path $D$ is uniquely determined by $\Set_X(D)$ and $\Set_Y(D)$ and furthermore
			$$\maj(D) = \sum_{i \in \Set_X(D)} i + \sum_{j \in \Set_Y(D)} j.$$
		\end{proposition}
		\begin{example}\label{ex2}
			Let $D = 01 \hspace{5pt} 001 \hspace{5pt} 011 \hspace{5pt} 01 \hspace{5pt} 01 \in \mathcal{D}_6$ be the Dyck path shown in Figure~\ref{fig3}. As indicated by the blanks, the descent set of $D$ is given by $\Des(D) = \{2,5,8,10 \}$ and its major index by $2+5+8+10 = 25$. The valleys of $D$ -- indicated in the picture by dots -- have coordinates $(1,1),(2,3),(4,4)$ and $(5,5)$ which gives
			$$ \Set_X(D) = \{1,2,4,5\}, \quad \Set_Y(D) = \{1,3,4,5\}.$$
		\end{example}
		There exist plenty of bijections between pattern-avoiding permutations and Dyck paths. To mention one, J.~Bandlow and K.~Killpatrick introduced an interesting bijection which sends the \Dfn{inversion number} $\inv(\sigma) := |\{i<j : \sigma_i > \sigma_j\}|$ of a $312$-avoiding permutation $\sigma$ to the area statistic of the corresponding Dyck path. We will discuss this and other bijections in more detail in Sections~\ref{sec5} and~\ref{sec6}, for definitions see \cite{bandlowkillpatrick}.
		
	\section{The bijection}\label{sec3}
	
		In this section, we construct a bijection $\Phi$ between $\Sn(231)$ and $\Dn$ such that
		$$\maj(\Phi(\sigma)) = \maj(\sigma)+\imaj(\sigma).$$
		\begin{lemma}\label{descentinverse}
			Let $\sigma \in \Sn(231)$, with $\Des(\sigma) = \{i_1,\dots,i_k\}$, $\Asc(\sigma) = \{j_1,\dots,j_{n-k}\}$. Then
			\begin{align*}
				\iDes(\sigma) &= \{\sigma_{i_1}-1,\dots,\sigma_{i_k}-1\}, \\
				\iAsc(\sigma) &= \big(\{\sigma_{j_1}-1,\dots,\sigma_{j_{n-k}}-1\}\setminus \{0\} \big) \cup \{n\}
			\end{align*}
		\end{lemma}
		\begin{proof}
			Let $i$ be a descent of $\sigma$. As $\sigma$ is $231$-avoiding, $[\sigma_i-1,\sigma_i]$ cannot be a subword of $\sigma$. In other words, $\hat{\sigma}(\sigma_i-1) > \hat{\sigma}(\sigma_i) = i$, and therefore $\sigma_i-1$ is a descent of $\hat{\sigma}$. On the other hand, let $i'$ be a descent of $\hat{\sigma}$, which is $312$-avoiding. The same argument as above yields $\sigma(\hat{\sigma}_{i'+1}) > \sigma(\hat{\sigma}_{i'+1}+1)$ or, equivalently, $\hat{\sigma}_{i'+1}$ is a descent of $\sigma$. This implies
			$$\iDes(\sigma) = \{\sigma_{i_1}-1,\dots,\sigma_{i_k}-1\}.$$
			As $\Asc(\sigma) = [n]\setminus\Des(\sigma)$, the statement about $\iAsc(\sigma)$ follows.
		\end{proof}
		
		Recall the definition of the involution $\rho$ on $\Sn$ and its properties as discussed in Section~\ref{sec2}. Together with Lemma~\ref{descentinverse}, this implies the following corollary.
		\begin{corollary}\label{cor1}
			Let $\sigma$ be $\tau$-avoiding for a given $\tau \in \{ 132, 231, 312, 213 \}$. Then 
			$$
				\des(\sigma) = \des(\hat{\sigma}), \quad \asc(\sigma) = \asc(\hat{\sigma}).
			$$
		\end{corollary}
		\begin{remark}
			For $\tau \in \{ 123, 321 \}$, the analogous statement of the previous corollary is false. E.g., $\sigma=[2,4,1,3]$ is $123$-avoiding and
			$$
				\Des(\sigma) = \{2\}, \quad \iDes(\sigma) = \{1,3\}.
			$$
		\end{remark}
		\begin{lemma}\label{lemma1}
			Let $j$ be an ascent of a $231$-avoiding permutation $\sigma$ and let $k > j$. Then $\sigma_k > \sigma_j$.
		\end{lemma}
		\begin{proof}
			This follows immediately from the fact that $\sigma$ is $231$-avoiding.
		\end{proof}
		For the remaining part of this section, set $\sigma$ to be a $231$-avoiding permutation. Our next goal is to show that $\sigma$ is uniquely determined by its ascent set and the image of its ascent set. Given the sets
		\begin{align*}
			\Asc(\sigma) = \{j_1,\dots,j_{n-k}=n\}, \quad \sigma(\Asc(\sigma)) = \{\ell_1,\dots,\ell_{n-k}\}
		\end{align*}
		in increasing order (i.e., $j_1 < \dots < j_{n-k}$ and $\ell_1 < \dots < \ell_{n-k}$). Then, by the previous lemma, $\sigma(j_1) < \dots < \sigma(j_{n-k})$, and therefore $\sigma(j_i) = \ell_i$.

		To determine $\sigma$ on its descent set, we determine $\sigma$ on its \Dfn{descent blocks}, i.e., the maximal sequences of consecutive descents. On any descent block $\{j_i+1,\dots,j_{i+1}-1\}$, $\sigma$ is decreasing and bounded from below by the image on the ascent following the sequence (recall that our definition of the ascent set implies that every descent block is followed by an ascent). In symbols, $\sigma(j_i+1) > \dots > \sigma(j_{i+1}-1) > \sigma(j_{i+1})$. Together with the property of being $231$-avoiding, this determines $\sigma$ on its descents from right to left: start with the rightmost descent, say $m$. Then $\sigma(m)$ is given by the smallest unassigned integer larger than $\sigma(m+1)$. Repeat this process by going through all descents from right to left as shown in Example~\ref{ex1}.

		By Lemma~\ref{descentinverse} and the above discussion we conclude the following fact, which is, in a slightly different context, originally due to A.~Reifegerste, see \cite[Proposition 2.4]{reifegerste} and the following discussion.
		\begin{corollary}\label{weneedit2}
			$\sigma$ is uniquely determined by $\Des(\sigma)$ and $\iDes(\sigma)$.
		\end{corollary}
		
		\begin{example}\label{ex1}
			Let $\sigma \in \mathcal{S}_6(231)$ such that
			$$
				\Des(\sigma) = \{1,2,4,5\}, \quad \iDes(\sigma) = \{1,3,4,5\}.
			$$
			This implies
			$$
				\Asc(\sigma) =\{3,6\}, \quad \iAsc(\sigma) = \{2,6\}.
			$$
			By Lemma~\ref{descentinverse}, we have
			$$
				\sigma(\Des(\sigma)) = \{2,4,5,6\}, \quad \sigma(\Asc(\sigma)) =  \{1,3\}.
			$$
			As described in the previous discussion, we first determine $\sigma$ on its ascent set,
			$$
				\sigma_3 = 1 < \sigma_6 = 3.
			$$
			Then, we determine $\sigma$ on its descents from right to left. The rightmost descent, $5$, is mapped to the smallest unassigned integer larger than $\sigma(6) = 3$, which is $4$; the next descent, $4$, is mapped to the smallest unassigned integer larger than $\sigma(5) = 4$, which is $5$; and so on. This gives
			$$
				\sigma_6 = 3 < \sigma_5 = 4 < \sigma_4 = 5, \quad \sigma_3 = 1 < \sigma_2 = 2 <	\sigma_1 = 6,
			$$
			and in total, we get $\sigma = [6,2,1,5,4,3]$.
		\end{example}
		Our next goal is to construct a bijection between $\Asc(\sigma)$ and $\iAsc(\sigma)$ such that the image of any ascent $j$ is less than or equal to $j$.
		
		Let $j$ be an ascent of $\sigma$. By $\tau(j)$, we denote  the size of the descent block immediately left of $j$, or, equivalently,
		$$\tau(j) := j-1-j',$$
		with $j'$ being the largest ascent such that $j'<j$ (respectively $0$ if $j$ is the first ascent). 
		\begin{lemma}\label{lemma2}
			Let $j$ be an ascent of $\sigma$. Then
			$$j \geq \sigma_j+\tau(j).$$
		\end{lemma}
		\begin{proof}
			Lemma~\ref{lemma1} implies $\sigma_k > \sigma_j$ for all $k$ such that $k>j$ and for all $k$ such that $j > k > j-\tau(j)$. Therefore, $n-j \leq n-\sigma(j)-\tau(j)$, which is equivalent to the statement.
		\end{proof}
		\begin{corollary}\label{inequality}
			Let $\Asc(\sigma) = \{j_1,\dots,j_{n-k}\}$. Then, for any $1 \leq \ell < n-k$, we have
			$$j_\ell \geq \sigma_{j_{\ell+1}}-1.$$
		\end{corollary}
		\begin{proof}
			By Lemma~\ref{lemma2}, $j_{\ell+1}$ is greater than or equal to $\sigma_{j_{\ell+1}} + \tau(j_{\ell+1})$, which, by definition, is equal to $\sigma_{j_{\ell+1}} + j_{\ell+1} - 1 - j_\ell$. This proves the corollary.
		\end{proof}
		By Lemma~\ref{descentinverse} and Corollary~\ref{inequality}, we can define a bijection between $\Asc(\sigma)$ and $\iAsc(\sigma)$, which has the desired property that the image of an ascent $j$ is less than or equal to $j$ in the following way:
		\begin{align*}
			j_\ell  &\mapsto \sigma_{j_{\ell+1}}-1 \quad \text{for} \quad 1 \leq \ell < n-k, \\
			j_{n-k} &\mapsto n
		\end{align*}
		This implies the following corollary concerning the descent sets of $\sigma$ and $\hat{\sigma}$:

		\begin{corollary}\label{weneedit}
			Let $\Des(\sigma) = \{i_1,\dots,i_k\}$ and let $\iDes(\sigma) = \{i'_1,\dots,i'_k\}$ such that $i_1 < \dots < i_k$ and $i'_1 < \dots < i'_k$. Then $i_\ell \leq i'_\ell$ for $1 \leq \ell \leq n$.
		\end{corollary}
		
		Now we are in the position to define the proposed bijection.

		\begin{definition}\label{bijectiondef}
			Define a map $\Phi$ from $\Sn(231)$ to $\Dn$ as follows: let $\sigma \in \Sn(231)$ with $\Des(\sigma) = \{i_1,\dots,i_k\}$ and $\iDes(\sigma) = \{i'_1,\dots,i'_k\}$, such that
			$$i_1 < \dots < i_k < i_{k+1} := n, \quad i'_1 < \dots < i'_k < i'_{k+1} := n.$$
			Then $\Phi(\sigma)$ is defined to be the Dyck path starting with $i'_1$ $0$'s, followed by $i_1$ $1$'s, followed by $i'_2-i'_1$ $0$'s, followed by $i_2-i_1$ $1$'s, followed by $i'_3-i'_2$ $0$'s, and so on, ending with $i_{k+1}-i_k$ $1$'s.
		\end{definition}
		
		\begin{theorem}\label{maintheorem}
			The map $\Phi$ defined in the previous definition is well-defined and bijective.
		\end{theorem}
		
		\begin{proof}
			By Proposition~\ref{cor1} and Corollary~\ref{weneedit}, $\Phi$ is well-defined and by Proposition~\ref{weneedit2} it is injective and therefore bijective.
		\end{proof}

		From the definition, we immediately obtain the following corollary.
		\begin{corollary}
			Let $\sigma \in \Sn(231)$. Then
			$$
				\Set_X(\Phi(\sigma)) = \Des(\sigma), \quad \Set_Y(\Phi(\sigma)) = \iDes(\sigma),
			$$
			 and, in particular,
			$$\maj(\Phi(\sigma)) = \maj(\sigma)+\imaj(\sigma).$$
		\end{corollary}

		\begin{example} \label{ex4}
			Let $\sigma$ be the permutation defined in Example~\ref{ex1}. As the descent set of $\sigma$ is $\{1,2,4,5\}$ and the descent set of its inverse is $\{1,3,4,5\}$, the coordinates of the descents of $\Phi(\sigma)$ are $(1,1),(2,3),(4,4)$ and $(5,5)$ 	and therefore $\Phi(\sigma) = 01 \hspace{5pt} 001 \hspace{5pt} 011 \hspace{5pt} 01 \hspace{5pt} 01$ is the Dyck path described in Example~\ref{ex2}. Furthermore,
			\begin{align*}
				\maj(\sigma)  &= 1+2+4+5 \hspace{5pt} = 12, \\
				\imaj(\sigma) &= 1+3+4+5 \hspace{5pt} = 13, \\
				\maj(\Phi(\sigma)) &= 2+5+8+10=25.
			\end{align*}
			In Figure~\ref{fig1} all $231$-avoiding permutations of the set $\{1,2,3,4\}$ and their associated Dyck paths are shown.
		\end{example}
		\begin{figure}
			\centering

			\setlength{\unitlength}{0.85pt}

			\begin{picture}(500,100)
			  \linethickness{.25\unitlength}
			  \color{grau}
			
				\multiput(0,0)(100,0){5}{
					\put(0 ,0 ){\line(0,1){80}}
					\put(20,20){\line(0,1){60}}
					\put(40,40){\line(0,1){40}}
					\put(60,60){\line(0,1){20}}
					
					\put( 0,20){\line(1,0){20}}
					\put( 0,40){\line(1,0){40}}
					\put( 0,60){\line(1,0){60}}
					\put( 0,80){\line(1,0){80}}
			
					\put(0,0){\line(1,1){80}}
				}
			
				\color{schwarz}
			  \linethickness{1.25\unitlength}
				% Erster pfad
				\put(0  ,0 ){\line(0,1){20.625}}
				\put(0  ,20){\line(0,1){20.625}}
				\put(0  ,40){\line(0,1){20.625}}
				\put(0  ,60){\line(0,1){20.625}}
					
				\put(0  ,80){\line(1,0){20.625}}
				\put(20 ,80){\line(1,0){20.625}}
				\put(40 ,80){\line(1,0){20.625}}
				\put(60 ,80){\line(1,0){20.625}}
			
				\put(0,90){\hbox{$[1,2,3,4] \mapsto$}}
			
				% Zweiter pfad
				\put(100,0 ){\line(0,1){20.625}}
				\put(100,20){\line(0,1){20.625}}
				\put(100,40){\line(0,1){20.625}}
				\put(120,60){\line(0,1){20.625}}
					
				\put(100,60){\line(1,0){20.625}}
				\put(120,80){\line(1,0){20.625}}
				\put(140,80){\line(1,0){20.625}}
				\put(160,80){\line(1,0){20.625}}
			
				\put(100,90){\hbox{$[4,1,2,3] \mapsto$}}
			
				% Dritter pfad
				\put(200,0 ){\line(0,1){20.625}}
				\put(200,20){\line(0,1){20.625}}
				\put(200,40){\line(0,1){20.625}}
				\put(240,60){\line(0,1){20.625}}
					
				\put(200,60){\line(1,0){20.625}}
				\put(220,60){\line(1,0){20.625}}
				\put(240,80){\line(1,0){20.625}}
				\put(260,80){\line(1,0){20.625}}
			
				\put(200,90){\hbox{$[1,4,2,3] \mapsto$}}
			
				% Vierter pfad
				\put(300,0 ){\line(0,1){20.625}}
				\put(300,20){\line(0,1){20.625}}
				\put(300,40){\line(0,1){20.625}}
				\put(360,60){\line(0,1){20.625}}
					
				\put(300,60){\line(1,0){20.625}}
				\put(320,60){\line(1,0){20.625}}
				\put(340,60){\line(1,0){20.625}}
				\put(360,80){\line(1,0){20.625}}
			
				\put(300,90){\hbox{$[1,2,4,3] \mapsto$}}
			
				% Fuenfter pfad
				\put(400,0 ){\line(0,1){20.625}}
				\put(400,20){\line(0,1){20.625}}
				\put(420,40){\line(0,1){20.625}}
				\put(420,60){\line(0,1){20.625}}
					
				\put(400,40){\line(1,0){20.625}}
				\put(420,80){\line(1,0){20.625}}
				\put(440,80){\line(1,0){20.625}}
				\put(460,80){\line(1,0){20.625}}
			
				\put(400,90){\hbox{$[3,1,2,4] \mapsto$}}
			
			\end{picture}

			\vspace{20pt}

			\begin{picture}(500,100)
			  \linethickness{.25\unitlength}
			  \color{grau}
			
				\multiput(0,0)(100,0){5}{
					\put(0 ,0 ){\line(0,1){80}}
					\put(20,20){\line(0,1){60}}
					\put(40,40){\line(0,1){40}}
					\put(60,60){\line(0,1){20}}
					
					\put( 0,20){\line(1,0){20}}
					\put( 0,40){\line(1,0){40}}
					\put( 0,60){\line(1,0){60}}
					\put( 0,80){\line(1,0){80}}
			
					\put(0,0){\line(1,1){80}}
				}
			
				\color{schwarz}
			  \linethickness{1.25\unitlength}
				% Erster pfad
				\put(0  ,0 ){\line(0,1){20.625}}
				\put(0  ,20){\line(0,1){20.625}}
				\put(20 ,40){\line(0,1){20.625}}
				\put(40 ,60){\line(0,1){20.625}}
					
				\put( 0,40){\line(1,0){20.625}}
				\put(20,60){\line(1,0){20.625}}
				\put(40,80){\line(1,0){20.625}}
				\put(60,80){\line(1,0){20.625}}
			
				\put(0,90){\hbox{$[4,3,1,2] \mapsto$}}
			
				% Zweiter pfad
				\put(100,0 ){\line(0,1){20.625}}
				\put(100,20){\line(0,1){20.625}}
				\put(120,40){\line(0,1){20.625}}
				\put(160,60){\line(0,1){20.625}}
					
				\put(100,40){\line(1,0){20.625}}
				\put(120,60){\line(1,0){20.625}}
				\put(140,60){\line(1,0){20.625}}
				\put(160,80){\line(1,0){20.625}}
			
				\put(100,90){\hbox{$[4,1,3,2] \mapsto$}}
			
				% Dritter pfad
				\put(200,0 ){\line(0,1){20.625}}
				\put(200,20){\line(0,1){20.625}}
				\put(240,40){\line(0,1){20.625}}
				\put(240,60){\line(0,1){20.625}}
					
				\put(200,40){\line(1,0){20.625}}
				\put(220,40){\line(1,0){20.625}}
				\put(240,80){\line(1,0){20.625}}
				\put(260,80){\line(1,0){20.625}}
			
				\put(200,90){\hbox{$[1,3,2,4] \mapsto$}}
			
				% Vierter pfad
				\put(300,0 ){\line(0,1){20.625}}
				\put(300,20){\line(0,1){20.625}}
				\put(340,40){\line(0,1){20.625}}
				\put(360,60){\line(0,1){20.625}}
					
				\put(300,40){\line(1,0){20.625}}
				\put(320,40){\line(1,0){20.625}}
				\put(340,60){\line(1,0){20.625}}
				\put(360,80){\line(1,0){20.625}}
			
				\put(300,90){\hbox{$[1,4,3,2] \mapsto$}}
			
				% Fuenfter pfad
				\put(400,0 ){\line(0,1){20.625}}
				\put(420,20){\line(0,1){20.625}}
				\put(420,40){\line(0,1){20.625}}
				\put(420,60){\line(0,1){20.625}}
					
				\put(400,20){\line(1,0){20.625}}
				\put(420,80){\line(1,0){20.625}}
				\put(440,80){\line(1,0){20.625}}
				\put(460,80){\line(1,0){20.625}}
			
				\put(400,90){\hbox{$[2,1,3,4] \mapsto$}}
			
			\end{picture}

			\vspace{20pt}

			\begin{picture}(400,100)
			  \linethickness{.25\unitlength}
			  \color{grau}
			
				\multiput(0,0)(100,0){4}{
					\put(0 ,0 ){\line(0,1){80}}
					\put(20,20){\line(0,1){60}}
					\put(40,40){\line(0,1){40}}
					\put(60,60){\line(0,1){20}}
					
					\put( 0,20){\line(1,0){20}}
					\put( 0,40){\line(1,0){40}}
					\put( 0,60){\line(1,0){60}}
					\put( 0,80){\line(1,0){80}}
			
					\put(0,0){\line(1,1){80}}
				}
			
				\color{schwarz}
			  \linethickness{1.25\unitlength}
				% Erster pfad
				\put(0  ,0 ){\line(0,1){20.625}}
				\put(20 ,20){\line(0,1){20.625}}
				\put(20 ,40){\line(0,1){20.625}}
				\put(40 ,60){\line(0,1){20.625}}
					
				\put( 0,20){\line(1,0){20.625}}
				\put(20,60){\line(1,0){20.625}}
				\put(40,80){\line(1,0){20.625}}
				\put(60,80){\line(1,0){20.625}}
			
				\put(0,90){\hbox{$[4,2,1,3] \mapsto$}}
			
				% Zweiter pfad
				\put(100,0 ){\line(0,1){20.625}}
				\put(120,20){\line(0,1){20.625}}
				\put(120,40){\line(0,1){20.625}}
				\put(160,60){\line(0,1){20.625}}
					
				\put(100,20){\line(1,0){20.625}}
				\put(120,60){\line(1,0){20.625}}
				\put(140,60){\line(1,0){20.625}}
				\put(160,80){\line(1,0){20.625}}
			
				\put(100,90){\hbox{$[2,1,4,3] \mapsto$}}
			
				% Dritter pfad
				\put(200,0 ){\line(0,1){20.625}}
				\put(220,20){\line(0,1){20.625}}
				\put(240,40){\line(0,1){20.625}}
				\put(240,60){\line(0,1){20.625}}
					
				\put(200,20){\line(1,0){20.625}}
				\put(220,40){\line(1,0){20.625}}
				\put(240,80){\line(1,0){20.625}}
				\put(260,80){\line(1,0){20.625}}
			
				\put(200,90){\hbox{$[3,2,1,4] \mapsto$}}
			
				% Vierter pfad
				\put(300,0 ){\line(0,1){20.625}}
				\put(320,20){\line(0,1){20.625}}
				\put(340,40){\line(0,1){20.625}}
				\put(360,60){\line(0,1){20.625}}
					
				\put(300,20){\line(1,0){20.625}}
				\put(320,40){\line(1,0){20.625}}
				\put(340,60){\line(1,0){20.625}}
				\put(360,80){\line(1,0){20.625}}
			
				\put(300,90){\hbox{$[4,3,2,1] \mapsto$}}
			
			\end{picture}

			\caption{The bijection between $\mathcal{S}_4(231)$ and $\mathcal{D}_4$.}
			\label{fig1}
		\end{figure}
	
	\section{A bistatistic on $231$-avoiding permutations}\label{sec4}
	
		Define the polynomial $\Aqt$ by
		$$\Aqt := \sum_{\sigma \in \Sn(231)}q^{\maj(\sigma)}t^{\binom{n}{2}-\imaj(\sigma)}.$$
		The first values of $\Aqt$ are given by
		\begin{align*}
			\mathcal{A}_1(q,t) &= 1, \\
			\mathcal{A}_2(q,t) &= q+t, \\
			\mathcal{A}_3(q,t) &= q^3+q^2t+qt^2+t^3+qt, \\
			\mathcal{A}_4(q,t) &= q^6+q^5t+q^4t^2+2q^3t^3+q^2t^4+qt^5+t^6+q^4t+q^3t^2+q^2t^3+qt^4+q^3t+qt^3.
		\end{align*}
		Obviously, $\Aqt$ reduces for $q=t=1$ to $\Cat_n$ and the bijection $\Phi$ defined in the previous section shows that $\Aqt$ can also be described in terms of Dyck paths. As considered by J.~F\"urlinger and J.~Hofbauer in \cite[Section 5]{fuerlingerhofbauer}, define two statistics $\maj_0$ and $\maj_1$ on a Dyck path $D$ by
		$$
			\maj_0(D) := \sum_{i \in \Des(D)}\big|\{j \leq i : D_j = 0 \}\big|, \quad \maj_1(D) := \sum_{i \in \Des(D)}\big|\{j \leq i : D_j = 1 \}\big|.
		$$
		\begin{corollary} Let $\Phi$ be the bijection defined in Theorem~\ref{maintheorem}, and let $\sigma$ be a $231$-avoiding permutation. Then
			$$
				\maj_1(\Phi(\sigma)) = \maj(\sigma), \quad \maj_0(\Phi(\sigma)) = \imaj(\sigma),
			$$
			and furthermore
			$$\Aqt = \sum_{D \in \Dn}q^{\maj_1(D)}t^{\binom{n}{2}-\maj_0(D)}.$$
		\end{corollary}
		Since $\maj_1(D)+\maj_0(D) = \maj(D)$, we get
		$$
			q^{\binom{n}{2}} \Aqq = \sum_{D \in \Dn}q^{\maj(D)},
		$$
		which is equal to MacMahon's $q$-Catalan number defined in the introduction, see \cite[Section 3]{fuerlingerhofbauer}. Using another identity in \cite{fuerlingerhofbauer}, we obtain the following generating function identity.
	
		\begin{theorem}\label{gfidentity}
			$\Aqt$ satisfies the generating function identity
			$$\sum_{n \geq 0}\frac{\Aqt \hspace{3pt} z^n}{(1+qz)\dots(1+q^{n+1}z)(1+tz)\dots(1+t^{n+1}z)} = 1.$$
		\end{theorem}
		\begin{proof}
			Setting $x=1, a=q^{-1}$ and $b=q$ in \cite[Theorem 5]{fuerlingerhofbauer} gives the above identity.
		\end{proof}
		The following corollary follows immediately from the fact that the generating function identity proved in Theorem~\ref{gfidentity} is symmetric in $q$ and $t$.
		\begin{corollary} \label{involution}
			The polynomial $\Aqt$ is symmetric in $q$ and $t$,
			$$\Aqt = \Atq.$$
		\end{corollary}
		Now, we want to provide a bijective proof of a refinement of this symmetry.
		\begin{theorem}\label{involution2}
			$$\sum_{\sigma \in \Sn(231)}a^{\des(\sigma)}q^{\maj(\sigma)}t^{\imaj(\sigma)} = \sum_{\sigma \in \Sn(312)}a^{n-1-\des(\sigma)}q^{\binom{n}{2} - \maj(\sigma)}t^{\binom{n}{2}-\imaj(\sigma)}.$$
		\end{theorem}
		\begin{proof}
			Let $\Psi$ be the involution on $\Sn(231)$ defined by the rule that a given $\sigma$ is mapped to the unique $\tau$ with $\Des(\tau):=[n-1] \setminus \iDes(\sigma)$ and $\iDes(\tau):=[n-1] \setminus \Des(\sigma)$ (in Section~\ref{sec3}, we proved the existence of such a $\tau$).
			This implies $\des(\sigma) = n-1-\des(\Psi(\sigma))$ and furthermore
			$$
				\maj(\sigma) = \binom{n}{2} - \imaj (\Psi(\sigma)), \quad \imaj(\sigma) = \binom{n}{2} - \maj(\Psi(\sigma)).
			$$
			Since mapping a permutation to its inverse interchanges $\Des(\sigma)$ and $\iDes(\sigma)$, the statement follows.
		\end{proof}
		\begin{remark}\label{psidef2}
			Equivalently, we could have defined the bijection $\Psi$ in the proof of Theorem~\ref{involution2} in terms of Dyck paths by the rule that a given Dyck path $D$ is mapped to the unique Dyck path $D'$ with
			$$
				\Set_X(D') = [n-1] \setminus \Set_Y(D), \quad \Set_Y(D') = [n-1] \setminus \Set_X(D).
			$$
		\end{remark}
	\section{Another description of $\Phi$}\label{sec5}
		The bijection $\Phi$ defined in Definition~\ref{bijectiondef} is closely connected to other bijections from pattern-avoiding permutations to Dyck paths. In \cite{krattenthaler2}, C. Krattenthaler constructed bijections from $132$- and from $123$-avoiding permutations to Dyck paths which were recently related to other bijections between pattern-avoiding permutations and Dyck paths by D. Callan in \cite{callan}. In this section, we express the bijection $\Phi$ in terms of Krattenthaler's bijection from $132$-avoiding permutations to Dyck paths, which we denote by $\kappa$.
		
		First, we recall the definition of $\kappa$ from \cite[Section 2]{krattenthaler2} and give an example. For any Dyck path $D$, the \Dfn{height of} $D$ \Dfn{at position} $i$ is the number of north steps minus the number of east steps until and including position $i$. Now, let $\sigma = [\sigma_1,\dots,\sigma_n]$ be a $132$-avoiding permutation and let $h_i$ denote the number of $j$'s larger than $i$ such that $\sigma_j > \sigma_i$. Read $\sigma$ from left to right and successively generate a Dyck path. When $\sigma_i$ is read, then in the path we adjoin as many north steps as necessary, followed by an east step from height $h_i+1$ to height $h_i$.
		\begin{example} \label{ex3}
			Let $\sigma = [3,4,5,1,2,6]$. Then $(h_1,\dots,h_6) = (3,2,1,2,1,0)$ which gives the heights of the east steps of the corresponding Dyck paths. Therefore, $\kappa(\sigma) = 0000111 \hspace{5pt} 00111$.
		\end{example}
		Our next goal is to express $\kappa$ in terms of the descent set of $\sigma$ and the descent set of $\hat{\sigma}$.
		\begin{lemma}
			Let $\sigma$ be any permutation. Then
			\begin{itemize}
				\item[(i)] $h_{i+1} < h_i$ if and only if $i \in \Asc(\sigma)$,
				\item[(ii)] $\sigma \in \Sn(132)$ if and only if $h_{i+1} \geq h_i - 1$ for all $1 \leq i < n$.
			\end{itemize}
		\end{lemma}
		\begin{proof}
			The first statement is obvious. as is the fact that $\sigma \in \Sn(132)$ implies $h_{i+1} \geq h_i - 1$ for all $1 \leq i < n$. For the reverse statement we use that $h_{i+1} \geq h_i - 1$ implies for an ascent $i$ that there exists no $k > i$ with $\sigma_i < \sigma_k < \sigma_{i+1}$.
		\end{proof}
		\begin{proposition}
			Let $\sigma$ be a $132$-avoiding permutation. Then $\kappa(\sigma)$ is the Dyck path $D$ given by
			$$
				\Set_X(D) = \Des(\sigma), \quad \Set_Y(D) = \{i+h_i : i \in \Des(\sigma)\}.
			$$
		\end{proposition}
		\begin{proof}
			The statement follows from the first part of the previous lemma and the definition of $\kappa$.
		\end{proof}		
		\begin{proposition}
			Let $\sigma$ be a $132$-avoiding permutation. Then
			$$\iDes(\sigma) = \big\{n-i-h_i : i \in \Des(\sigma) \big\}.$$
		\end{proposition}
		\begin{proof}
			Let $i$ be a descent of $\sigma$. The number of $\sigma_j$'s right of $\sigma_i$ that are smaller than $\sigma_i$ is equal to $n-i-h_i$. Since $\sigma$ is $132$-avoiding, we know that the set of those $\sigma_j$'s is equal to $\{1,\dots,n-i-h_i\}$. This implies that $\hat{\sigma}_{n-i-h_i} > \hat{\sigma}_{n-i-h_i+1}$, and therefore $n-i-h_i$ is a descent of $\hat{\sigma}$.
			
			Now let $i$ be a descent of $\sigma$ and let $j > i$. The previous lemma implies that $i+h_i \neq j+h_j$, because otherwise $\{i,i+1,\dots,j-1\} \subseteq \Asc(\sigma)$, a contradiction to $i \in \Des(\sigma)$. The proposition follows from Corollary~\ref{cor1}.
		\end{proof}
		These two propositions characterize $\kappa$ in terms of the descent set of $\sigma$ and the descent set of $\hat{\sigma}$: let $\sigma \in \Sn(132)$. Then
		$$
			\Set_X(\kappa(\sigma)) = \Des(\sigma), \quad \Set_Y(\kappa(\sigma)) = \{n-j : j \in \iDes(\sigma)\}.
		$$
		\begin{example}
			Continuing the previous example, we have that $\Des(\sigma)$ and $\iDes(\sigma)$ are given by $\{3\}$ and by $\{2\}$, respectively. On the other hand,
			$$
				\Set_X(\kappa(\sigma)) = \{3\}  =  \Des(\sigma), \quad \Set_Y(\kappa(\sigma)) = \{4\} = \{n-j : j \in \iDes(\sigma)\}.
			$$
		\end{example}
		Let $D$ be a Dyck path with valleys $\{(i_1,j_1),\dots,(i_k,j_k)\}$. The last map we need is the involution $\Psi$ as described in Remark~\ref{psidef2} followed by the involution on Dyck paths sending $D$ to the Dyck path with valleys $\{(n-j_k,n-i_k),\dots,(n-j_1,n-i_1)\}$. We denote this composition by $\psi$.
		
		Now we can describe the relation of $\Phi$ and $\kappa$ using the involutions $\rho$ and $\psi$:
		\begin{theorem}
			$\kappa = \psi \circ \Phi \circ \rho$.
		\end{theorem}
		\begin{proof}
			Let $\sigma$ be a $132$-avoiding permutation. Then
			\begin{align*}
				\Set_X(\Phi \circ \rho (\sigma)) &= \Des(\rho(\sigma))  = [n-1] \setminus \{n-i : i \in \Des(\sigma) \}, \\
				\Set_Y(\Phi \circ \rho (\sigma)) &= \iDes(\rho(\sigma)) = [n-1] \setminus \iDes(\sigma).
			\end{align*}
			On the other hand,
			\begin{align*}
				\Set_X(\psi(D)) &= \{n-i : i \in [n-1] \setminus \Set_X(D) \}, \\
				\Set_Y(\psi(D)) &= \{n-i : i \in [n-1] \setminus \Set_Y(D) \},
			\end{align*}
			for a Dyck path $D$. So, in total, we have
			\begin{align*}
				\Set_X(\psi \circ \Phi \circ \rho (\sigma)) &= \Des(\sigma), \\
				\Set_Y(\psi \circ \Phi \circ \rho (\sigma)) &= \{n-j : j \in \iDes(\sigma)\}.
			\end{align*}
			Together with the above discussion, this implies the theorem.
		\end{proof}
		\begin{example}
			Set $\sigma := [6,2,1,5,4,3]$. As we have seen in Example~\ref{ex4}, $\Phi(\sigma) = 01 \hspace{3pt} 001 \hspace{5pt} 011 \hspace{5pt} 01 \hspace{5pt} 01$, and therefore $\psi(\Phi(\sigma)) = 0000111 \hspace{5pt} 00111$. As shown in Example~\ref{ex3}, this equals $\kappa(\rho(\sigma))$.
		\end{example}
		Next, we want to use the involutions $\Psi, \rho$, and $\sigma \mapsto \hat{\sigma}$ (we denote the latter by $\i$) to prove an analogue of Theorem~\ref{involution2} for $\Sn(132)$ and $\Sn(213)$.
		\begin{theorem}
			$$\sum_{\sigma \in \Sn(132)}a^{\des(\sigma)}q^{\maj(\sigma)}t^{\imaj(\sigma)} = \sum_{\sigma \in \Sn(213)}a^{n-1-\des(\sigma)}q^{\binom{n}{2} - \maj(\sigma)}t^{\binom{n}{2}-\imaj(\sigma)}.$$
		\end{theorem}
		\begin{proof}
			The bijection between $\Sn(132)$ and $\Sn(213)$ defined by $\rho \circ \i \circ \Psi \circ \rho$, where $\Psi$ is meant as described in the proof of Theorem~\ref{involution2}, sends the tristatistic $(\des,\maj,\imaj)$ to the tristatistic $(n-1-\des,\binom{n}{2} - \maj, \binom{n}{2} - \imaj)$:
			\begin{align*}
				(\des,\maj,\imaj) &\xmapsto{\rho} (n-1-\des,\binom{n}{2} - n \des+\maj,\binom{n}{2}-\imaj) \\
													&\xmapsto{\psi} (\des,\imaj, n \des - \maj) \\
													&\xmapsto{\i} 	(\des,n \des - \maj,\imaj) \\
													&\xmapsto{\rho} (n-1-\des,\binom{n}{2} - \maj,\binom{n}{2} - \imaj).
			\end{align*}
		\end{proof}
		\begin{corollary}
			$$\sum_{\sigma \in \Sn(132)}q^{\maj(\sigma)}t^{\imaj(\sigma)} = \sum_{\sigma \in \Sn(213)}q^{\binom{n}{2} - \maj(\sigma)}t^{\binom{n}{2}-\imaj(\sigma)}.$$
		\end{corollary}

		The analogous statement is true if summing over $\Sn(123)$ and over $\Sn(321)$ respectively. This can be easily deduced from the following theorem, which was conjectured in a former version of this article\footnote{The author is very grateful to Matteo Silimbani for providing a proof of Theorem~\ref{th:silimbani}.}.

		\begin{theorem}\label{th:silimbani}
			There exists an involution on $\Sn(321)$ which leaves the descent set $\Des$ invariant and maps the descent set of the inverse, $\iDes$, to $\{n-j : j \in \iDes\}$.
		\end{theorem}
		\begin{proof}
			By the Robinson-Schensted correspondence, permutations in $\Sn$ are in bijection with pairs $(P,Q)$ of standard Young tableaux of shape $\lambda \vdash n$, see e.g. \cite[Chapter 19]{foatahan}. It is well-known that $\Des(\sigma) = \Des(Q)$ and $\iDes(\sigma) = \Des(P)$ and that the length of the first column in the shape of $P$ and $Q$ equals the length of a longest decreasing chain in the associated permutation. In particular, we have that $321$-avoiding permutations are in bijection with pairs of standard Young tableaux with at most 2 rows. The last ingredient we need is that the \emph{promotion operator} $P \mapsto P^J$ is an involution, leaves the shape invariant and has the property that $\Des(P^J) = \{n-j : n \in \Des(P)\}$, see the discussion after Theorem 19.3 in \cite{foatahan}. Therefore, the involution $\mathbf{j}$ on $\Sn(321)$ defined by
			$$\sigma \mapsto (P,Q) \mapsto (P^J,Q) \mapsto \mathbf{j}(\sigma)$$
			 has the desired property.
		\end{proof}
		
		In the remaining part of this section, we describe the connection between $\Phi$ and the bijection by Bandlow and Killpatrick between $\Sn(132)$ and $\Dn$ mentioned at the end of Section~\ref{sec2}. This shows how the area statistic on Dyck paths can be represented in $\Sn(231)$ via $\Phi$. We denote the bijection Bandlow  and Killpatrick defined in \cite{bandlowkillpatrick} by $\beta$. Together with M.~Fulmek, C.~Krattenthaler already described the connection between $\beta$ and $\kappa$ in a comment on Bandlow and Killpatrick's article, \cite{fulmekkrattenthaler}. By the definitions of $\Phi$, $\kappa$ and $\beta$, we have
		$$\beta = \Psi \circ \Phi \circ \i.$$
		As $\beta$ sends the inversion number on $\Sn(132)$ to the area statistic on $\Dn$ and $\i$ leaves the inversion number invariant, we conclude for any $231$-avoiding permutation $\sigma$,
		$$\operatorname{inv}(\sigma) = \area(\Psi \circ \Phi(\sigma)).$$
	
	\section{Connections to $q,t$-Catalan numbers}\label{sec6}
	
		The $q,t$-Catalan numbers $\Cat_n(q,t)$ can be defined in terms of the two statistics \emph{area} and \emph{bounce} on Dyck paths, which were first considered by J.~F\"urlinger and J.~Hofbauer in \cite{fuerlingerhofbauer} and by J.~Haglund in \cite{haglund2} respectively. For definitions and for further information see e.g. \cite{haglund}.
		\begin{definition}
			Define the $q,t$\Dfn{-Catalan numbers} $\Cat_n(q,t)$ as
			$$\Cat_n(q,t) := \sum_{D \in \Dn}{q^{\operatorname{area}(D)}t^{\operatorname{bounce}(D)}}.$$
		\end{definition}
		In \cite{garsiahaiman}, A.~Garsia and M.~Haiman proved that $\Cat_n(q,t) = \Cat_n(t,q)$. But it is still an open problem to find a bijective proof of this fact. Furthermore, they proved that $q^{\binom{n}{2}} \Cat_n(q,q^{-1})$ is equal to MacMahon's $q$-Catalan numbers, and therefore
			$$\Cat_n(q,q^{-1}) = \Aqq.$$
		Simple computations show that, for $n \geq 4$, $\Cat_n(q,t) \neq \Aqt$, but the computations lead to the following guess. It was verified for $n \leq 10$.
		\begin{guess} \label{gue:Cat-A}
			We have
			$$\Cat_n(q,t) = \sum_{D \in \Dn}q^{\maj_1(D)-k_D}t^{\binom{n}{2}-\maj_0(D)-k_D},$$
			for appropriate non-negative integers $k_D$.
		\end{guess}
		\begin{example}
			Let $n=4$. Then
			$$\Cat_4(q,t) = q^6+q^5t+q^4t^2+q^3t^3+q^2t^4+qt^5+t^6+q^4t+q^3t^2+q^2t^3+qt^4+q^3t+q^2t^2+qt^3,$$
			and therefore
			$$\mathcal{A}_4(q,t) - \Cat_4(q,t) = q^3t^3-q^2t^2 = q^2t^2(qt-1).$$
			This shows that
			$$\Cat_4(q,t) = \sum_{D \in \mathcal{D}_4}q^{\maj_1(D)-k_D}t^{\binom{n}{2}-\maj_0(D)-k_D},$$
			where either $k_D = \Big\{ \begin{array}{rl} 1,& D = 01 \hspace{5pt} 01 \hspace{5pt} 0011 \\ 0,& \text{otherwise} \end{array}$, or $k_D = \Big\{ \begin{array}{rl} 1,& D = 000111 \hspace{5pt} 01 \\ 0,& \text{otherwise} \end{array}$.
		\end{example}
		\begin{remark}
			For $n=4$, the involution defined in Corollary~\ref{involution} interchanges the paths $01 \hspace{5pt} 01 \hspace{5pt} 0011$ and $000111 \hspace{5pt} 01$. Therefore even if one could determine the $k_D$'s in Guess~\ref{gue:Cat-A}, the proposed involution would fail to prove the open problem of finding a bijective proof of the symmetrie property of $\Cat_n(q,t)$.
		\end{remark}
%		Finally, we want to describe a refinement of $\Cat_n(q,t)$ in terms of $231$-avoiding permutations. For that, we also need the fact that
%		$$\des(\sigma) = n-1-\des(\Psi \circ \Phi(\sigma)).$$
%		In \cite[Theorem 3.15]{haglund}, J.~Haglund describes a bijection on Dyck paths which sends the area statistic to the bounce statistic. We denote this bijection by $\alpha$. From the definition it is very easy to deduce that for any Dyck path $D \in \Dn$,
%		$$\des(D) = n-1-\des(\alpha(D)).$$
%		Using this bijection as well as the bijections $\Phi$ and $\Psi$, we can describe
%		$$\Cat_n(a,q,t) := \sum_{D \in \Dn}{a^{\des(D)}q^{\operatorname{area}(D)}t^{\operatorname{bounce}(D)}}$$
%		in terms of $231$-avoiding permutations. For the sake of more readability, set
%		$$\star := \Phi^{-1} \circ \Psi^{-1}.$$
%		\begin{corollary}
%			We have
%			\begin{align*}
%		\Cat_n(a,q,t) &= \sum_{D \in \Dn}{a^{\des(\star \circ \alpha^{-1}(D))}q^{\operatorname{inv}(\star(D))}t^{\operatorname{inv}(\star \circ \alpha^{-1}(D))}} \\
%									&= \sum_{\sigma \in \Sn(231)}{a^{\des(\sigma)}q^{\operatorname{inv}(\star \circ \alpha \circ \star^{-1}(\sigma))}t^{\operatorname{inv}(\sigma)}}.
%			\end{align*}
%		\end{corollary}
%		\begin{remark}
%			Of course, this construction just translates Haglund's bijection on Dyck paths to $231$-avoiding permutations. It does not give any further hint how the symmetry problem could be solved.
%		\end{remark}

		\providecommand{\bysame}{\leavevmode\hbox to3em{\hrulefill}\thinspace}
\providecommand{\MR}{\relax\ifhmode\unskip\space\fi MR }
\providecommand{\MRhref}[2]{%
  \href{http://www.ams.org/mathscinet-getitem?mr=#1}{#2}
}
\providecommand{\href}[2]{#2}

%		\bibliographystyle{amsplain}
%		\bibliography{../bibliography}
	
\end{document}